\documentclass[11pt]{article}
\usepackage{hyperref}
\usepackage{amssymb,amsmath, amsthm}
\usepackage{graphicx}
\usepackage{graphics}
\usepackage{psfrag}

\usepackage[all]{xy}

\newtheorem{theorem}{Theorem}[section]

\theoremstyle{definition}

\newcommand{\Rb}{\mathbb{R}}

\begin{document}
 
\title{Numerical Solution of the Neural Field Equation in the Two-Dimensional Case}
\author{Pedro M. Lima and Evelyn Buckwar\\
Institute of Stochastics\\
Johannes Kepler University\\
Altenbergerstr. 69 , 4040 Linz, Austria}
\maketitle

\begin{abstract}
We are concerned with the numerical solution of a class integro-differential equations, known as Neural Field Equations, which  describe the large-scale dynamics of spatially structured networks of neurons. These equations have many applications in Neuroscience and Robotics.
We describe a numerical method for the approximation of solutions in the two-dimensional case, including a time-dependent delay in the integrand function. Compared with known algorithms for this type of equation we propose a scheme with higher accuracy in the time discretisation. Since computational efficiency is a key issue in this type of  calculations, we use  a new method for reducing the complexity of the algorithm.
The convergence issues are discussed in detail and a number of numerical examples is presented, which illustrate the performance of the method.
\end{abstract}

\section{Introduction}

 In recent years significant progress has been made in understanding the brain electrodynamics using mathematical techniques. Neural field models represent the large-scale dynamics of spatially structured networks of neurons in terms of nonlinear integro-differential equations. Such models are becoming increasingly important for the interpretation of experimental data, including those obtained from EEG, fMRI and  optical imaging \cite{Coombes}. These equations also play an important role in Cognitive Robotics, since the architecture of autonomous robots, able to interact with other agents in solving a mutual task, is strongly inspired by the processing principles and the neuronal circuitry in the primate brain (see \cite{Erlhagen}).
All this explains why Neural Field Equations {\em in several dimensions} are a very actual and important subject of research.

Moreover, simulations play a fundamental role in studying brain dynamics in Computational Neuroscience, and to understand diseases such as Parkinson, as well as the effect of treatments, such as in Deep Brain Stimulations (DBS) or Transcranial Magnetic Stimulations (TMS). Thus, the availability of efficient, fast, reliable numerical methods is an important ingredient for improving the effectiveness of techniques such as DBS or TMS in many therapeutic applications.
 Integro-differential equations in several spacial dimensions are quite a challenge for numerical simulation, because the standard approaches require a very high computational effort. Maybe this explains why  there exist very few studies concerning the numerical analysis of NFEs. Some important papers which inspired our work are \cite{Faye,Potthast, Hutt2014}; as other relevant references we can cite  \cite{BL,CV} .  Overall the lack of efficient  algorithms represents a serious drawback for the use of NFEs in practical applications. This is the main motivation for the present work.

We are concerned with the numerical solution of the following integro-differential
equation:
\begin{equation} \label{1}
c \frac{\partial}{\partial t} V(\bar{x},t) =
I(\bar{x},t) - V(\bar{x},t) + \int_{\Omega} K(|\bar{x}-\bar{y}|) S(V(\bar{y},t))d \bar{y},
\end{equation}
\[ t\in[0,T], \bar{x} \in \Omega \subset \Rb^2,\]
where the unknown $V(\bar{x},t)$ is a continuous function 
$V: \Omega \times [0,T] \rightarrow \Rb$, $I$, $K$ and $S$ are given functions; $c$
is a constant. In this article, by $|\bar{x}-\bar{y}|$ we mean $\|\bar{x}-\bar{y}\|_2$.
We search for a solution $V$ of this equation which satisfies
the initial condition
\begin{equation} \label{1a}
V(\bar{x},0)= V_0(\bar{x}), \qquad \bar{x}  \in \Omega.
\end{equation} 
Along with equation \eqref{1} we will also consider
\begin{equation} \label{2}
c \frac{\partial}{\partial t} V(\bar{x},t) =
I(\bar{x},t) - V(\bar{x},t)+ \int_{\Omega} K(|\bar{x}-\bar{y}|) S(V(y,t-\tau(\bar{x},\bar{y}))d\bar{y},
\end{equation}
\[  t\in[0,T],  \quad  \bar{x}  \in \Omega \subset \Rb^2,\]
where $\tau(\bar{x},\bar{y} )>0$ is a delay, depending on the spatial variables.
In the latter case, the initial condition has the form
\begin{equation} \label{2a}
V(\bar{x},t)= V_0(\bar{x},t), \qquad \bar{x}  \in \Omega, \quad t \in [-\tau_{max},0],
\end{equation} 
where $\tau_{max}=\max_{\bar{x},\bar{y} \in \Omega} \tau (\bar{x},\bar{y})$.
By integrating both sides of \eqref{2} with respect to time on $[0,T]$, we obtain
the Volterra-Fredholm integral equation:
\begin{equation} \label{3}
c  V(\bar{x},t) = V_0(\bar{x})+ \int_0^t \left((I(\bar{x},s) - V(\bar{x},s)+ \int_{\Omega} K(|\bar{x}-\bar{y}|) S(V(\bar{y},s-\tau(\bar{x},\bar{y}))d\bar{y}\right) ds,  
\end{equation}
\[  t\in[0,T], \bar{x}  \in \Omega \subset \Rb^2.\]
The existence and uniqueness of a solution of equation \eqref{1} in the case $\Omega= \Rb^m$, $m\ge 2$, was proved
in \cite{Potthast}, both in  the case of a smooth and discontinuous function $S$.
An analytical study of equation \eqref{2} was carried out in \cite{Faye},
where the authors have addressed the problems of existence, uniqueness and stability of solutions.
They define  the  Banach space $ \tilde{C}=C([l, F)$, were $l$ is some real interval, containing  0, and $F$ denotes the Banach space
$L^2(\Omega, \Rb)$,with the norm
$$
\| \Psi \|_{F}= \sqrt{\int_{\Omega} \Psi^2(r) dr}, \quad \forall \Psi \in  F.
$$
Hence the norm in $ \tilde{C}$ is defined by
$$
\| \Phi \|_{\tilde{C}} = \sup_{t \in l}\sqrt{ \int_{\Omega} \Phi^2(r,t) dr}, \quad \forall \Phi \in  \tilde{C}.
$$
They have proved that if $K \in L^2(\Omega^2,\Rb)$, 
$I \in \tilde{C}([-\tau_{max},\infty[,F )$
and $\tau \in C (\bar{\Omega}^2,\Rb_+) $, then for every $V_0 \in \tilde{C} ([-\tau_{max}, 0])$ equation \eqref{2} has a unique solution 
 $V \in \tilde{C}^1([0, \infty[,F) \cap \tilde{C} ([-\tau_{max}, \infty[,F)$.
Subsequently in this paper we will use the same notations and definitions of norms. Moreover, as the authors of \cite{Faye},
we will assume that  $S$ and its derivative $S'$ are positive and bounded.
%%%%%%%%%%%%%%%%%%%%%%%%%%%%%%%%%%%%%%%%%%%%%%%%%%%%%%%%%%%%%%%%%%%%%%%%%%%%%%%%%
When solving numerically equations of the forms \eqref{1} and \eqref{2},
they are often reduced to the form \eqref{3}; therefore we begin by 
discussing literature on computational methods for Volterra-Fredholm equations.
 Starting with the one-dimensional case, without delay,
 Brunner has analysed the convergence of collocation methods \cite{Brunner}, while Kauthen
 has proposed continuous time collocation methods \cite{Kauthen}. In \cite{Han} an asymptotic error expansion for the Nystr\"{o}m method was proposed, which enabled the use
 of extrapolation algorithms to accelerate the convergence of the method.
 Another approach was developped by Z. Jackiewicz and co-authors \cite{Jacki1}, \cite{Jacki2}, who have applied Gaussian quadrature rules and interpolation
 to approximate the solution of integro-differential equations
 modelling neural networks, which are similar to equation \eqref{1}.
 
 In all the above mentioned papers the authors were concerned  with the one-dimensional case.
 In the two-dimensional case, the required computational effort to solve
 equations  \eqref{1} and \eqref{3}  grows very fast as
 the discretization step is reduced, and therefore special attention has to be paid 
 to the creation of effective methods. This can be achieved by means of 
 low-rank methods, as those discussed in \cite{Xie}, when the kernel is approximated
 by polynomial interpolation,  which enables a significant reduction of the dimensions of the matrices. In \cite{Cardone}, the authors use an iterative method
 to solve linear systems of equations which takes into account the special form of the matrix to introduce parallel computation.
%%%%%%%%%%%%%%%%%%%
Concerning equation \eqref{2},  besides the existence and stability of
 solution, numerical approximations were obtained in \cite{Faye}.
The computational method applies  quadrature rule in space to reduce the problem
to a system of delay differential equations, which is then solved
by a standard algorithm for this kind of equations. A more efficient approach was recently
proposed in \cite{Hutt2010} \cite{Hutt2014}, where the authors introduce a new approach to deal with
the convolution kernel of the equation and use Fast Fourier Transforms to reduce 
significantly the computational effort required by numerical integration.

The above mentioned equations are known as Neural Field Equations (NFE) and  have played an important role in mathematical
neuroscience for a long time. Equation \eqref{1} was introduced first by Wilson
and Cowan \cite{WC}, and then by Amari \cite{Amari}, to describe excitatory
and inhibitory interactions in populations of neurons.
While in other mathematical models of neuronal interactions  the function
$V$ (membrane potential) depends only on time, in the case of NFE it is
a function of time and space.
The function $I$ represents external sources of excitation and
$S$ describes the dependence between the firing rate of the neurons and their membrane
potential. It can be either a smooth function (typically of sigmoidal type)
or a Heaviside function. The kernel function $K(|\bar{x}-\bar{y}|)$ gives the connectivity
between neurons in positions  $\bar{x}$ and $\bar{y}$. By writing the arguments of the function
in this form we mean that we consider the connectivity homogeneous, that
is, it depends only on the distance between neurons, and not on their specific
location.

According to many authors, realistic models of neural fields must take
into account that the propagation speed of neuronal interactions is finite,
which leads to NFE with delays of the form \eqref{2}.

 In the present paper we propose a new numerical approach to the Neural Field Equation,
in the forms \eqref{1} and \eqref{3}. One remarkable feature of our method
is that we use use a implicit second order scheme for the time discretisation, which improves 
 its accuracy and stability, when compared with the available algorithms. Moreover, to reduce
the computational complexity of our method we use an interpolation procedure which allows a drastic
 reduction of matrix dimensions, without a significant loss of accuracy. 
This improves the efficiency of the algorithm.

The  outline of the method is as follows. In Sec. 2, we describe the numerical algorithm, 
for equation \eqref{1}; its stability, convergence and complexity are analysed. In the same
section, we introduce an algorithm for equation  \eqref{2}.
In Sec. 3  a set of numerical examples is presented, which illustrate both cases, and the numerical results are
discussed. We finish with some conclusions in Sec. 4.
%%%%%%%%%%%%%%%%%%%%%%%%%%%%%%%%%%%%%%%%%%%%%%%%%%%%%%%%%%%%%%%%%%%%%%%%%%%%%%%%%%%%%%%
\section{Numerical Method}
\subsection{Neural Field Equation without delay}
\subsubsection{Time Discretization}
We begin by rewriting equation \eqref{1}  in the form
\begin{equation} \label{21}
c \frac{\partial}{\partial t} V(\bar{x},t) =
I(\bar{x},t) - V(\bar{x},t)+ \kappa (V(\bar{x},t))
\end{equation}
\[ t\in[0,T], \bar{x}  \in \Omega \subset \Rb^2,\]
where $\kappa $ denotes the nonlinear integral operator defined by
\begin{equation} \label{22} 
\kappa (V(\bar{x},t))=\int_{\Omega} K(|\bar{x}-\bar{y}|) S(V(\bar{y},t))d\bar{y}.
\end{equation}
We shall first deal with the time discretization in equation \eqref{21}, therefore
we introduce the stepsize $h_t>0$ and define
\[ t_i = i h_t, \quad i=0,...,M, \quad T=h_t M. \]
Moreover, let
\[ V_i(\bar{x})= V(t_i,\bar{x}), \quad \forall x \in \Omega, \quad i=0,...,M.\]
We shall approximate the partial derivative in time by the backward difference
\begin{equation} \label{23} 
\frac{\partial}{\partial t} V(\bar{x},t_i) \approx \frac{3 V_i(\bar{x})-4 V_{i-1}(\bar{x})+V_{i-2}(\bar{x})}{2 h_t},
\end{equation}
which gives a discretization error of the order $O(h_t^2)$, for sufficiently smooth $V$. By substituting
\eqref{23} into \eqref{21} we obtain the implicit scheme
\begin{equation} \label{24}
c \frac{3 U_i-4 U_{i-1}+U_{i-2}}{ 2h_t}  =
I_i - U_i+ \kappa (U_i), \quad i=2,...,M,
\end{equation}
where  $U_i$  approximates the solution of \eqref{21}. 

To start this scheme we need to know $U_0$, which is defined by the initial condition $V_0$,
and $U_1$, which can be obtained by a one-step method, for example, the explicit Euler method.

It is important to analyse the {\em  stability} of the numerical scheme \eqref{24}. 
If we denote 
\[  e_i= \|V_i -U_i\| \]
(the discretization error at time $t_i$), we want to analyse the conditions
under which, if $I_i=0$,  we have $ e_i < e_{i-1} $, if $h_t$ is sufficiently small. 
 We can look at the scheme \eqref{24} as a multistep method for the solution of a system of differential
equations (the  approximate solution $U_i$ , at time $t_i$, is obtained from $U_{i-1}$
and $U_{i-2}$). Therefore, we  should begin by studying the zero-stability. If we multply both sides of
 \eqref{24} by $2 h_t$ and then ignore the terms containing $h_t$, we have  an equation of the form
\[ 3 U_i - 4U_{i-1} + U_{i-2} =0, \]
to which corresponds the characteristic equation
\[ 3  \mu^2 - 4 \mu +1=0.\]
The roots of this equation are 
\[ \mu_{1,2} = \frac{2 \pm 1 }{3};\]
it is easily seen that $|\mu_{1,2}| \le 1$ and $\mu_1=1$ is not a multiple root. Therefore, the scheme \eqref{24}
is zero-stable. This means that  we  have 
$e_i \le e_{i-1}$, if $h_t$ is sufficiently small.

The above result can be summarized in the form of the following theorem:
\begin{theorem}
The scheme \eqref{24} is zero-stable.
\end{theorem}

Our next step is to investigate under which conditions equation \eqref{24} has a unique solution, so that each step
of the iterative process is well defined.
With this purpose we write this equation in the form

\begin{equation} \label{24a}
  U_i (\bar{x})-\frac{ 1}{1+ \frac{3c}{2 h_t}} \kappa(U_i)= f_i(\bar{x}),  \qquad \bar{x} \in \Omega
\end{equation}
where
\begin{equation} \label{24b}
f_i(\bar{x}) = \left( 1+ \frac{2h_t}{3c}\right)^{-1} \left( I_i+\frac{c}{h_t} 2 U_{i-1}(\bar{x}) -\frac{c}{2h_t} U_{i-2}(\bar{x})\right), \qquad \bar{x} \in \Omega
\end{equation}
Equation \eqref{24a} - \eqref{24b} is a nonlinear Fredholm integral equation of the second kind and we
will analyse its solvability using standard results of functional analysis.
%so to prove its solvability we will apply similar arguments to those used in \cite{Davis}, p. 424, where 
%this type of equations is analysed in the one dimensional case.

In order to apply the Banach fixed point theorem, we define the iterative process:
\begin{equation} \label{25a}
  U_i	^{(\nu)} (\bar{x})=\lambda \kappa(U_i^{(\nu-1)})+ f_i(\bar{x}) =G(U_i^{\nu-1}) ,  \qquad \bar{x} \in \Omega, \quad n=1,2,...
\end{equation}
were
\begin{equation} \label{lam}
\lambda= \frac{ 1}{1+ \frac{3c}{2h_t}} = \frac{ 2 h_t}{2 h_t+ 3c}.
\end{equation}

If the function $G$ is contractive in a certain closed set  $X \subset F$, such that  $G(X) \subset X$,
then by the Banach fixed point theorem equation \eqref{24a} has a unique solution in $X$
and  the sequence $U_i^{(n)},$ defined by \eqref{25a}, converges to this solution in the norm of $F$,
for any initial guess $U_i^{(0)} \in X$. In our case, the solution is by construction the iterate $U_i$,
so it should be close to $U_{i-1}$ and $U_{i-2}$. Therefore it makes sense to assume that $X$ is  a certain set containing 
 $U_{i-1}$ and $U_{i-2}$ and to choose  $U_i^{(0)}=U_{i-1}$.

To prove that $G$ is contractive in $X$ we need to show that for a certain constant $L$, $L<1$, we
have
\begin{equation} \label{26a}
\| G(V)-G(U) \|_F \le L \| V-U \|_F,  \qquad \forall U,V \in X. 
\end{equation}
By definition
\begin{equation} \label{27a}
(\| G(V)-G(U) \|_F)^2 = \lambda^2 \int_{\Omega} |K( \bar{x}-\bar{y})|^2 |S(V)-S(U)|^2 d \bar{y}; 
\end{equation}
Using the mean value theorem for integrals, we  get
\begin{equation} \label{28a}
(\| G(V)-G(U) \|_F)^2 \le \lambda^2  |\Omega| (\|K\|_{L^2(\Omega^2)})^2 \max_{U,V \in X}(\|S(V)-S(U)\|_F)^2 ,
\end{equation}
where $|\Omega|$ denotes the area of $\Omega$. Since we have assumed that $S$ has a bounded continuous  derivative  in $\Rb$, we can write
\begin{equation} \label{29a}
\left(\| G(V)-G(U) \|_F \right)^2  \le |\Omega|\lambda^2  \left(\|K\|_{L^2(\Omega^2)} \right)^2  S_{max}^2  \left(\|V-U\|_F\right)^2 ,
\end{equation}
where
\[ 
S_{max} = \max_{x \in \Rb} |S'(x)| .
\]

Finally, we rewrite \eqref{29a} in the form
\begin{equation} \label{30a}
\| G(V)-G(U) \|_F \le \lambda  \sqrt{|\Omega|}\|K\|_{L^2(\Omega^2)} S_{max} \|V-U\|_F
\end{equation}
and we conclude that \eqref{26a} holds with
\begin{equation} \label{31a}
L= \lambda  \sqrt{|\Omega|}\|K\|_{L^2(\Omega^2)} S_{max} .
\end{equation}
Recall that
\[
\lambda= \frac{ 2 h_t}{2 h_t+ 3c} <  \frac{ 2 h_t}{ 3c}. 
\]
Then, in order to satisfy $L<1$ it is sufficient to require that
\begin{equation} \label{32a}
\frac{ 2 h_t}{ 3c} \sqrt{|\Omega|}\|K\|_{L^2(\Omega^2)} S_{max} < 1
\end{equation}
or equivalently
\begin{equation} \label{33a}
h_t < \frac{3c} {2\sqrt{|\Omega|}\|K\|_{L^2(\Omega^2)} S_{max} }.
\end{equation}
From \eqref{33a} we conclude that $G$ will be contractive in a certain set $X \subset F$ 
if we take $h_t$ sufficiently small.
Moreover, if $h_t$ is sufficiently small,  we can choose a set $X$ such that 
 $G(X) \subset X$ and $\{U_{i-1}, U_{i-2}\} \subset X$. 
Therefore the  iterative process  \eqref{25a} with  $U_i^{(0)}=U_{i-1}$ will converge to the
solution of \eqref{24a}.

The above construction not only shows that the equation (\ref{24a}) has a unique solution in a certain set $X $, but
it also suggests that the iterative process  \eqref{25a}, starting with  $U_i^{(0)}=U_{i-1}$  can be effectively used to approximate this solution.
Actually, the convergence rate of this process depends on the constant  $L$, which as follows from
\eqref{31a} is approximately proportional to $h_t$. In other words, the convergence of the process will be
faster and faster as $h_t$ tends to zero.

The above result can be formulated in the form of the following theorem.

\begin{theorem}
For each $i=2, 3...$ , if $h_t$ satisfies \eqref{33a} the nonlinear equation \eqref{24a} has a unique solution $U_i \in X$, where $X \subset F$ is a certain closed set containing $U_{i-1}$ and $U_{i-2}$.  Moreover,  the iterative process
\eqref{25a}  with $U_i^{(0)}=U_{i-1}$converges to this solution.
\end{theorem}

%%%%%%%%%%%%%%%%%%%%%%%%%%%%%%%%%%%%%%%%%%%%%%%%%%%%%%%%%%%%%%%%%%%%%%%%%%%%%%%%%%%%%%%%%%%%%%%%%
\subsubsection{Space Discretization}
Since the equation \eqref{24a} in general cannot be solved analytically,
we need a computational method to compute a numerical approximation of its solution.
By other words, we need a space discretization, which will be the subject of this subsection.  

 For the sake of simplicity, assume that $\Omega$ is a rectangle: $\Omega=[-1,1] \times [-1,1]$.
 We now introduce a uniform
grid of points $(x_i,x_j)$, such that \\ $x_i=-1+i h$,   $i=0,...,n$, where  $h$ is the discretization step
in space.  In each subinterval $[x_i,x_{i+1}]$   we introduce $k$
Gaussian nodes: $x_{i,s}=x_i + \frac h 2 ( 1+ \xi_s )$, $i=0,1,\dots n-1$, where $\xi_s$ are the roots of the $k$-th degree
Legendre polynomial, $s=1,...,k$.
We shall denote $\Omega_h$ the set of all grid points  $(x_{is}, x_{jt})$, $i,j=0,...,n-1$,$s,t=1,...,k$.
  A Gaussian quadrature formula to evaluate the integral 
$\int_{\Omega} f(u,v) du dv $ will have the form
  \begin{equation} \label{28}
Q(f) =  \sum_{i=0}^{n-1} \sum_{j=0}^{n-1} \sum_{s=1}^{k} \sum_{t=1}^{k} \tilde{w}_s \tilde{w}_t   
f(x_{is}, x_{jt}),
\end{equation}
 with $\tilde{w}_s =\frac h 2  w_s$, where $w_s$ are the standard weights of a Gaussian quadrature formula with $k$ nodes on $[-1,1]$, $s=1,...,k$.
As it is well-known, a quadrature formula of this type has degree $2k-1$ and therefore, assuming that $f$
has at least $2k$ continuous derivatives on $\Omega$, 
the integration error of \eqref{28} is of the order of $h^{2k}$.   
Note that the total number of nodes in the space discretization is $k^2 n^2$ . 

When we introduce the quadrature formula \eqref{28} to compute $\kappa(U)$ we define a finite-dimensional approximation
of the operator $\kappa$. Let us denote $U^h$ a vector with $N^2$ entries, where $N=n\,k$, such that
\[
(U^h)_{is,jt} \approx U(x_{is},x_{jt});
\]
then the finite-dimensional approximation of $\kappa(U)$ may be given by
 \begin{equation} \label{28b}
(\kappa^h(U^h)) _{mu,lv}=  \sum_{i=0}^{n_1} \sum_{j=0}^{n_2} \sum_{s=1}^{k} \sum_{t=1}^{k} \tilde{w}_s \tilde{w}_t   
K(\|(x_{mu},x_{lv})-(y_{is},y_{jt})\|_2) S((U^h)_{is,jt}). 
\end{equation}
By replacing $\kappa$  with $\kappa^h$ in equation \eqref{24a} we obtain the following system of nonlinear equations:
 \begin{equation} \label{28c}
  U^h -\frac{ 1}{1+ \frac{2h_t}{3c}} \kappa^h(U^h)= f^h,
\end{equation}
where $\kappa^h (U^h)$ is defined by \eqref{28b} and 
\[
(f^h)_{is,jt} =f(x_{is},x_{jt}),
\]
with $f$ defined by \eqref{24b}; in \eqref{28c}, for the sake of simplicity, we have omitted the index $i$ of $U_i^h$.
Note that for the the computation of $f^h$ we have to evaluate the iterates $U_{i-1}$ 
and $U_{i-2}$ at all the points of $\Omega_h$. We denote the vectors resulting from this evaluation by
$U_{i-1}^h$ and $U_{i-2}^h$, respectively.
We conclude that at each time step of our numerical scheme we must solve \eqref{28c}, which is
a system of $N^2$ nonlinear equations. From this discretization some questions arise:
\begin{enumerate}
\item Is the system \eqref{28c} solvable?
\item  Does the solution $U^h$  of this system converge in some sense to $U_i$, as $h \rightarrow 0$?
\item  How can we estimate the error $E^h_i = \| U^h-U_i\|$ ?
\end{enumerate}
We can investigate the solvability of  \eqref{28c} in the same way as we have studied the Fredholm
integral equation \eqref{24a}. More precisely, we can introduce the iterative procedure
\begin{equation} \label{27a}
  U^{h,(\nu)} =\lambda \kappa^h (U^{h,(\nu-1)})+ f^h =G^h(U^{h,(\nu-1)}) ,
\end{equation} 
$m=1,2, \dots$. 
As in the case of the Fredholm integral equation, the convergence of the iterative procedure \eqref{27a} 
depends on the contractivity of $G^h$. Using the same arguments as in subsection 2.1.1, 
we obtain the following inequality, which is a finite-dimensional analogue of \eqref{30a}:
\begin{equation} \label{30b}
\| G^h(V)-G^h(U) \|_2 \le \lambda  K_{max} S'_{max}\|V-U\|_2 \sum_{i=0}^{n} \sum_{j=0}^{n} \sum_{s=1}^{k} \sum_{t=1}^{k} \tilde{w}_s \tilde{w}_t  ,
\end{equation}
for $U,V \in X^h  \subset \Rb^{N^2}$, where
\begin{equation} \label{Kh}
K_{max} =\max_{(x_{mu},x_{lv}),(y_{is},y_{jt}) \in \Omega_h} |K(\|(x_{mu},x_{lv})-(y_{is},y_{jt})\|_2)|.
\end{equation}
Since, by the construction of the Gaussian quadrature, 
\[
\sum_{i=0}^{n} \sum_{j=0}^{n} \sum_{s=1}^{k} \sum_{t=1}^{k} \tilde{w}_s \tilde{w}_t =1,
\]
we get that $G^h$ is Lipschitzian in $X^h$ with the Lipschitz constant
\begin{equation} \label{lh}
L_1= \lambda K_{\max} S_{\max},
\end{equation}
where $\lambda$ is defined by \eqref{lam}.
Finally we conclude that $G_h$ is contractive if 
\begin{equation} \label{33b}
h_t < \frac{3c} {2 K_{\max} S_{\max} }.
\end{equation}
% %  theorem 2.3
\begin{theorem}
For each $i=2, 3...$ , if  $S$ and $K$ are such that $S_{max}$ and $K_{max}$ exist and  $h_t$ satisfies \eqref{33b}, then the nonlinear equation \eqref{28c} has a unique solution $U^h \in X^h$, where $X^h \subset \Rb^{N^2}$ is a certain closed set containing $U_{i-1}^h$ and $U_{i-2}^h$.  Moreover,  the iterative process
\eqref{27a}  with $U^{h,(0)}=U_{i-1}^h$ converges to this solution.
\end{theorem}
Since we know that the equation \eqref{28c} is solvable under certain conditions and we even know
an iterative scheme to obtain its solution, we should now address the problem of the convergence of the
solution $U^h$ to $U_i$, as $h \rightarrow 0$.
  
With this purpose, we write equation \eqref{24a} at each grid point of $\Omega_h$:
\begin{equation} \label{34a}
  U_i (x_{is,jt})-\frac{ 1}{1+ \frac{2h_t}{3c}} \kappa(U_i(x_{is,jt}))= f_i(x_{is,jt}),  \qquad x_{is,jt} \in \Omega^h.
\end{equation}
Subtracting from this equation \eqref{28c}, we obtain
\begin{equation} \label{35a}
  U_i (x_{is,jt})-  U^{h}_{is,jt}=\lambda ( \kappa(U_i(x_{is,jt})) - \kappa^h( U^{h}_{is,jt})) ,
	\qquad i,j=0,...,n-1, \quad s,t=1,...,k.
\end{equation} 
We note that 
\begin{equation} \label{36a}
\kappa(U_i(x_{is,jt})) - \kappa^h( U^{h}_{is,jt}) =
 \kappa(U_i(x_{is,jt})) -\kappa^h( U_i (x_{is,jt}))+  (\kappa^h( U_i(x_{is,jt}) -\kappa^h( U^{h}_{is,jt})).
\end{equation}
Substituting \eqref{36a} into \eqref{35a} yields
\begin{equation} \label{37a}
  U_i (x_{is,jt})-  U^{h}_{is,jt} =
	\lambda ( \kappa(U_i(x_{is,jt})) -\kappa^h( U_i({is,jt}))+  \kappa^h ( U_i({is,jt})) -\kappa^h( U^{h}_{is,jt})),
\end{equation} 
$\qquad i,j=0,...,n-1, \quad s,t=1,...,k.$
From \eqref{37a}, we obtain 
\begin{equation} \label{37b}
  \|U_i -  U^{h}\|_{\infty} \le
	\lambda ( \|  \kappa(U_i) -\kappa^h( U_i))\|_{\infty}+  \|\kappa^h( U_i) -\kappa^h( U^{h})\|_{\infty}).
\end{equation} 
The second term in the last sum can be expressed in terms of $ \|U_i -  U^{h}\| $, if we take into account
that $\kappa^h$ is a Lipschitzian function, that is, there exists a certain constant $L_1$  such that
\begin{equation} \label{38a}
\|\kappa^h( U_i) -\kappa^h( U^{h})\|_{\infty} \le  L_1 \|U_i -  U^{h}\|_{\infty}, \forall U_i, U^h \in X^h,
\end{equation} 
where $L_1$ is defined by \eqref{lh}. 
Hence we may rewrite \eqref{37a} in the form
\begin{equation} \label{38a}
  \|U_i -  U^{h}\|_{\infty} (1 - \lambda L_1) \le
	\lambda   \|  \kappa(U_i) -\kappa^h( U_i)\|_{\infty}.
\end{equation} 
Recall that $\lambda$ is approximately proportional to $h_t$ and therefore we have $\lambda L_1 <1$
if we choose $h_t$ sufficiently small.
 
Finally, in order to evaluate $\|  \kappa(U_i) -\kappa^h( U_i))\|_{\infty}$, we must remember that $\kappa^h$ results
from the approximation of the integral $\kappa$ by a Gaussian quadrature rule, with stepsize $h$
and $k$ nodes at each subinterval. Therefore, if $U_i$, $K$ and $S$ are sufficiently smooth, there exists a certain $M>0$, which does not depend on $h$, such that
\begin{equation} \label{39a}
\|  \kappa(U_i) -\kappa^h( U_i))\|_{\infty} \le M h^{2k}.
\end{equation}
Finally, from \eqref{38a} and \eqref{39a} we conclude that  there exists such a constant $\tilde{M}$ that
\begin{equation} \label{40a}
  \|U_i -  U^{h}\|_{\infty} \le  \tilde{M} h^{2k}.
	\end{equation} 
The above results lead to the following theorem.
%%%%%%%%%%%%%%%%% theorem 2.4
\begin{theorem}
Under the conditions of Theorem 2.3,  the unique solution $U^h$ of \eqref{28c} converges to the solution $U_i$ of \eqref{24a}
as $h \rightarrow 0$. Moreover, if  $U_i$, $K$ and $S$ are sufficiently smooth, the following error estimate holds
\[  \|U_i -  U^{h}\|_{\infty} =O(h^{2k}) ,  \qquad  \makebox{as} \, h \, \rightarrow 0.\]
\end{theorem}
%%%%%%%%%%%%%%%%%%%%%%%%%%%%%%%%%%%%%%%%%%%%%%%%%%%%%%%%%%
\subsubsection{ Computational Implementation}
The above numerical algorithm for the approximate solution of the neural field equation
in the two-dimensional case was implemented by means of a MATLAB code.

The code has the following structure. After introducing the input data (step size in time and in space, 
initial condition $U_0$, error tolerance for the inner cycle, required number of steps in time) , there is an outer 
cycle that computes each vector $U^h$ , given $U_{i-1}^h$ and $U_{i-2}^h$, according to the multistep
method \eqref{24}. In order to initialize this cycle, besides $U_0$, we need $U_1^h$, which is
obtained by the explicit Euler method. More precisely, we compute
\begin{equation} \label{40}
  U_1^h= U_0+ \frac{h_t}{c} (I_0-U_0+\kappa^h (U_0)).
	\end{equation} 
We recall that at each step in time we must solve the nonlinear system of equations \eqref{24a}, which as suggested
above is obtained by means of the fixed point method, that is, we iterate the scheme \eqref{27a}, until the iterates
satisfy 
\[
 \|U^{h,(n)}-  U^{h,(n-1)}\|_{\infty} < \epsilon,
\]
for some given $\epsilon$. This is the inner cycle of our scheme. Typically, in all the examples we have
computed the number of iterations in the inner cycle is not very high (3-4, in general), confirming that
the fixed point method is an efficient way of solving the system \eqref{24a}.
To start the inner cycle we use an initial guess which is obtained from $U_{i-1}^h$ using again the Euler method:
\begin{equation} \label{40}
  U^{h,(0)} = U_{i-1}+ \frac{h_t}{c} (I_i-U_{i-1}^h+\kappa^h (U_{i-1}^h)).
\end{equation} 
Note that at each step of the inner cycle it is necessary to compute the function $\kappa_h$ at all the grid points.
From the computational point of view, this means that we must evaluate $N^2$ times a quadrature rule of
the form \eqref{28} (with $N^2$ nodes). Of course, this requires a high computational effort and the greatest
part of the computing time of our algorithm is spent in this process. 
Therefore, we pay special attention to reducing the computational cost at this stage.
In order to improve the efficiency of the numerical method, we apply the following technique, proposed 
in \cite{Xie} for the solution of two-dimensional Fredholm equations.
Assuming that the function $V$ is sufficiently smooth, we  can approximate it by an interpolating 
polynomial of a certain degree. As it is known from  the theory of approximation, the best approximation of a smooth function by an interpolating polynomial of degree $m$ is obtained if the interpolating points are the roots of the Chebyshev polynomial of degree $m$:
 \begin{equation} \label{29}
p_i^m= cos \left(\frac{(2i-1)\pi}{m} \right), \qquad i=1,...,m
\end{equation}
Our approach for reducing the matrices rank in our method consists in replacing the  solution
$V_i$ by its interpolating polynomial at the Chebyshev nodes in $\Omega $. If  $V_i$ is sufficiently
smooth, this produces a very small error and yields a very significant reduction of computational cost.
Actually, when computing the vector $\tilde{\kappa} (U_i)$ (see formula \eqref{22}) we have only
to compute $m^2$ components, one for each Chebyshev node on $[-1,1] \times [-1,1]$. Choosing $m$ much smaller than $n$, we thus obtain a significant computational advantage.

The procedure at each iteration is as follows. We compute the matrix $M$ such that 
\[ M_{i,j}= Q(V( p_i^m,p_j^m,t))  , \qquad  i=1,...,m,  j=1,...,m, \]
where $Q$ is the approximation of the integral $\kappa$, obtained by means of the quadrature
\eqref{28}, $p_i^m$ are the Chebyshev nodes, defined by \eqref{29}.

Then we have to perform the matrix multiplication
\begin{equation} \label{30}
\Lambda =C M C^T , 
\end{equation} 
where $C$ is the matrix defined by
\[ C_{ij} = c_{i-1} (p_j^m), \qquad i=1,...,m, \quad j=1,...,m; \]
here $c_k$ represents the scaled Chebyshev polynomial of degree $k$,
\[ c_k (x) = \delta_k \cos(k \,\makebox{arcos}(x) ), \quad  k=0,1,.. \]
with $ \delta_0= 1/\sqrt{n}$,  $ \delta_k= \sqrt{2} \delta_0 $, $k=1,...,m-1$.
The matrix $\Lambda$ contains the coefficients of the interpolating
polynomial of the solution (expanded in terms of scaled Chebyshev polynomials).
Finally,  in order to obtain the interpolated values of the solution at the Gaussian nodes,
we have to compute
\begin{equation} \label{31}
 T =P^T  \Lambda P, 
\end{equation}
where $P$ is the transformation matrix, given by
\[ P_{ij} = c_{i-1} (x_{(j)} ) ,  \qquad  i=1,...,m , \quad j =1,..., N. \]
Here $x_{(j)}$ represents each Gaussian node: $x_{(j)} = x_{i,s} $, if $ j=ik+s$.
Finally, the vector $U_i$ for the next time step (of size $N^2$ ) is obtained by copying $T$, 
row by row (note that $T$ is a matrix of dimension $N \times N$).
%%%%%%%%%%%%%%%%%%%%%%%%%%%%%%%%%%%%%%%%%%%%%%%%%%%%%%%%%%%%%%%%%%%%%%
\subsubsection{Complexity Analysis}
As remarked before, it is important to analyse the complexity of the computations, since
the computational effort can be signifficantly reduced by the application of adequate techniques.
In the previous section, we have described  an algorithm for computing each iterate
of the fixed point method, which requires $m^2$ applications of the quadrature formula \eqref{28}.
Since this quadrature implies $N^2$ evaluations of the integrand function,
we have a total of $m^2 N^2$ function evaluations.
Note that if no polynomial interpolation would be applied, $N^4$ evaluations of the integrand function would be required at each iteration. It is easy to conclude that the number of arithmetic operations
required to apply the quadrature is also proportional to $m^2N^2$.

Then, according to the described algorithm, we must perform the matrix multiplication
\eqref{30}. Since the involved matrices have dimension $m \times m$, the total
number of arithmetic operations is $O(m^3)$ . Since, by construction, $m <<N$, the
complexity of this part of the computations is much less than the previous one.
%This complexity can be further reduced if we remark that the matrix product 
%\eqref{30} corresponds to the application
%of a two-dimensional Discrete Cosinus Transform (DCT) to the matrix $M$.Then
%we can apply a fast  algorithm to perform this transformation (analagous to the algorithm
%of the Fast Fourier Transform) and the complexity can be reduced to $O(m^2 \log_2(m))$.

Finally, we have the matrix product $\eqref{31}$. Here the transformation matrix $P$
has dimensions  $m \times N$, as the resulting matrix $T$ has dimensions  $N \times N$
The resulting complexity is therefore  $O(m N^2)$.

In conclusion, the number of evaluations of the integrand function in each iterate
of the fixed point method is $N^2 m^2$ and the complexity of each iteration
is $O( N^2 m^2 )$. Note that the number of iterations of the fixed point method 
at each time step is typically $2-4$.
%%%%%%%%%%%%%%%%%%%%%%%%%%%%%%%%%%%%%%%%%%%%%%%%%%%%%%%%%%%%%%%%%%%%%%%%%%%%%%%%%%
\subsubsection{Error Analysis}
We start by analysing the error resulting from the time discretization.
Assuming that the partial derivatives $ \frac{ \partial ^i V(x,t)}{\partial^i t}$, $i=1,2,3$, are continuous
 on a certain domain $\Omega \times [0,T]$, the local discretization error
of the approximation, given by \eqref{23}, has the order of $O(h_t^2)$.

Concerning the space discretization, the error has two components: one resulting from the application
of the discretization scheme \eqref{28c} , and the other resulting from the polynomial interpolation.
In both cases, the order  of the approximation depends on the smoothness of the solution.
 Therefore we must choose the degree of the Gaussian quadrature  according to the
smoothness of the mentioned functions.

The first component was analysed in Sec. 2.1.2. According to Theorem 2.4,
if the functions $S$ , $K$ and $U_i$ satisfy certain smoothness conditions, the 
discretization scheme \eqref{28c} has order $2k$ in space, where $k$ is the number of Gaussian nodes at each subinterval. 

To analyse the interpolation error, we refer to Lemma 3 in \cite{Xie}.
According to this Lemma, if the partial derivatives $ \frac{ \partial ^i f(y_1,y_2)}{\partial^i y_j}$
of  a certain function $f$ are continuous, with $j=1,2$, $i=1,2,..,s$ then
 \begin{equation} \label{33}
 \|f -C_m f \| =O(m^{-s} \log^2 m) , 
\end{equation}
where $C_m f$ represents the interpolating $m$-th degree polynomial of $f$ in the Chebyshev nodes.

In order to obtain an optimal precision of the method, we require that the components
of the error, given by  \eqref{40a}  and \eqref{33}, are of the same order.  Hence, to ensure 
that we can choose  $m \ll N$, the degree of smoothness $s$ of the solution
should be significantly higher than $2k$ (otherwise we will not obtain 
a significant reduction of the matrices rank). For example, suppose that
$m=N^{1/2}=h^{-1/2}$, then for the interpolation error
we have  $m^{-s} \log^2 m =h^{s/2} \log^2 (h^{-1/2})$=$-\frac 1 2 h^{s/2} \log^2 h$ .
Comparing with \eqref{40a}, we conclude that for optimal precision
we must have $s/2 \approx 2k$.

Hence  the described method is specially suitable
in the case of a smooth solution , so that one can take advantage of the
small quadrature error and strong rank reduction.

Summarizing, we have so far shown that the numerical scheme has local
discretization order $O(h_t^2) +O(h^{2k})$. Let us now analyse the global error
 \[  E_{i,j}=V( \bar{x_j}, t_i) - (U^h_i)_j. \]
For $i=0$, we have $ E_{0,j}=0$,$j=1,...,N^2$. For $i=1$, we know that $U_1^h$ is computed according to \eqref{40}.
According to the results of section 2.1.2 about space discretization, and since the Euler method has local discretization error
of order 2, we have
\begin{equation} \label{global0}
 \| E_{1,j}\|=\max_{ \bar{x_j}\in \Omega_h} |V( \bar{x_j}, t_1) - (U^h_1)_j|=  O(h_t^2) +O(h^{2k}).
\end{equation} 
Concerning the subsequent steps in time, $i=2,3,...$ from \eqref{28c} and \eqref{24a} we can conclude that
\begin{equation} \label{global1}
 E_{i,j}-\frac 4 3 E_{i-1,j} +\frac 1 3 E_{i-2,j} = \frac{2 h_t}{3c} \left(E_{i,j} +\kappa(V_i)_j - \kappa^h(U^h_i)_j \right).
\end{equation}
As in Section 2.1.2, we can write 
\begin{equation} \label{global2}
\kappa(V_i)_j - \kappa^h(U^h_i)_j=\kappa(V_i)_j - \kappa^h(V_i)_j +\kappa^h(V_i)_j - \kappa^h(U^h_i)_j
\end{equation} 
and therefore
\begin{equation} \label{global3}
\|\kappa(V_i)_j - \kappa^h(U^h_i)_j\| \le  L_1\| E_{i,j}\| +O(h^{2k}),
\end{equation} 
where $L_1$ is given by \eqref{lh}. Substituting \eqref{global3} into \eqref{global1} we obtain
\begin{equation} \label{global4}
\| E_{i,j}\| \le \frac 4 3 \|E_{i-1,j}\| +\frac 1 3 \|E_{i-2,j}\|  +\frac{2 h_t}{3c} \left( \|E_{i,j}\| +
 L_1\| E_{i,j}\| +O(h^{2k}) \right),  \qquad i=2,3,...
\end{equation}
which is equivalent to 
\begin{equation} \label{global5}
\| E_{i,j}\|(1- \frac{2 h_t}{3c}(1+ L_1 ) )\le \frac 4 3 \|E_{i-1,j}\| +\frac 1 3 \|E_{i-2,j}\|+O(h^{2k}), \qquad i=2,3,...
\end{equation}
provided that
\begin{equation} \label{global6}
\frac{2 h_t}{3c}(1+ L_1) < 1.
\end{equation}
In particular, for $i=2$, we get
\begin{equation} \label{global8}
\| E_{2,j}\| \le \left(1- \frac{2 h_t}{3c}(1+ L) \right)^{-1} \frac 4 3 \|E_{1,j}\| +\frac 1 3 \|E_0\|+O(h^{2k}) ;
\end{equation}
according to \eqref{global0} we conclude that
\begin{equation} \label{global9}
\| E_{2,j}\| \le (1- \frac{2 h_t}{3c}(1+ L) )^{-1} ( O(h^{2k})+ O(h_t^{2})) +O(h^{2k}) = O(h_t^{2}) + O(h^{2k}) .
\end{equation}
Due to the complexity of our numerical method, it was not possible to obtain a closed expression of
the global error, for an arbitrary value of $i$. 
However, according to Theorem 2.1, the multistep scheme is zero-stable, which means that it will be stable for a sufficiently 
small value of $h_t$. The condition \eqref{global6} actually shows us how small $h_t$ must be in order to achieve stability.
If this condition is satisfied,  we expect that the scheme will be stable and the method will
be convergent, with the convergence order   $O(h_t^{2}) + O(h^{2k}) $  (the same as for the local discretization error). 

We remark that if we used an explicit  method, like the Euler method, the stability condition on $h_t$
would be much more restrictive. In other words, the fact that we use an implicit method allows us to use larger step size in time,
which makes the method more efficient.

The numerical results presented  in Sec. 3  are in agreement with the error analysis and confirm the expected convergence orders, for a set of different cases.  
%%%%%%%%%%%%%%%%%%%%%%%%%%%%%%%%%%%%%%%%%%%%%%%%%%%%%%%%%%%%%%%%%%%%%%%%%%%%%
\subsection{Delay Equation}
We now focus our attention on equation \eqref{2}, where the argument of the solution inside the integral
has a delay $\tau(\bar{x},\bar{y})$. This delay takes into account the fact that the propagation speed of signals between 
neurons is finite and therefore the post-synaptic potential generated at location $\bar{x}$ in instant $t$ by action potentials 
arriving from connected neurons at location $\bar{y}$ actually depends on the potential of these neurons at instant 
$t-\tau(\bar{x},\bar{y})$, where $\tau(\bar{x},\bar{y})$ is the time taken by the signal to come from $\bar{y}$ to
$\bar{x}$. Since we assume that the propagation speed $v$ is constant and uniform in space, we have
\[ \tau(\bar{x},\bar{y})=\frac{|\bar{y}-\bar{x}|}{v}.
\]
Hence, the delay integro-differential equation that we must solve has the form
\begin{equation} \label{delay}
c \frac{\partial}{\partial t} V(\bar{x},t) =
I(\bar{x},t) - V(\bar{x},t)+ \int_{\Omega} K(|\bar{x}-\bar{y}|) S(V(\bar{y},t-\frac{|\bar{y}-\bar{x}|}{v})) d\bar{y}.
\end{equation}
Note that in this case the initial conditions satisfied by the solution of our problem  have the form \eqref{2a},
where 
\[ \tau_{max} =\max_{\bar{x},\bar{y}\in \Omega} \frac{|\bar{y}-\bar{x}|} {v} .\]

The numerical algorithm used to solve equation \eqref{delay} is essentially the same as described in the previous
sections. The main difference results from the fact that when computing the integral on the right-hand side of \eqref{delay}
at instant $t_i$ we must use not only the approximate solution at instants $t_{i-1}$ and $t_{i-2}$, but at all
instants $t_{i-k}$, $k=1,...,k_{max}$, where $k_{max}$  is the integer part of $\tau_{max}/h_t$.
Note also that the argument $t_i -\frac{|\bar{y}-\bar{x}|}{v}$ may not be a multiple of $h_t$.
In general let $j$  and $\delta_t$ be the integer and the fractional part of  $\frac{|\bar{y}-\bar{x}|}{v h_t} $.
In this case, we have
\[ t_{i-j-1} \le  t_i -\frac{|\bar{y}-\bar{x}|}{v}  \le t_{i-j} 
\]
and
\[ h_t \delta_t= \left(t_i- \frac{|\bar{y}-\bar{x}|}{v} \right)-t_{i-j-1}.
\]
The needed value of the solution $V(\bar{y},t_i-\frac{|\bar{y}-\bar{x}|}{v})$ is then approximated by linear interpolation:
\begin{equation} \label{inter}
  V\left(\bar{y},t_i-\frac{|\bar{y}-\bar{x}|}{v}\right) \approx \delta_t U_{i-j} + (1-\delta_t ) U_{i-j-1}. 
\end{equation}
Note that the error analysis that we have carried out in the previous subsection may not apply to the delay equation.
What we can say in this case, assuming that $V$ is a smooth function of $t$, is that 
the error introduced each time we use the approximation formula \eqref{inter} has the order of $O(h_t^2)$ (the same
as the error resulting from the time discretization). However, the overall effect of this error in the computations requires a more 
detailed analysis, which is left for a future work.

Concerning complexity, in the case of the delay equation, each time we compute the integrand function, we
must compute the delay $\tau(\bar{x},\bar{y})$. As discussed above, this delay is obtained dividing  the distance $\|\bar{y}-\bar{x}\|_2$
by $v$. Since this distance is also required to compute  the kernel connectivity $K(\|\bar{y}-\bar{x}\|_2)$, for an effective computation
this quantity should be evaluated only once and then kept in memory.

\section{Numerical Results}
\subsection{Neural Field Equation without delay}
Here we present the results of some numerical tests we have carried out, in order to 
check the convergence properties of the described method (in the case where no delay is considered).
 Our main purpose is to test experimentally the convergence of the method and measure the error; therefore we have chosen
some cases where the exact solution is known and do not arise from applications. However the form of the
connectivity kernels and firing rate functions in these examples are close to the ones of neuroscience problems.
We first check the convergence order in time. With this purpose, we consider the following example.

{ \bf Example 1.} In this example, 
\[ K(|\bar{x}-\bar{y}|)=K(x_1,x_2,y_1,y_2) = \exp \left( -\lambda(x_1-y_1)^2 -\lambda(x_2-y_2)^2 \right), \]
where $\lambda \in \Rb^+$;
$S (x)= \tanh (\sigma x) $,  $\sigma \in \Rb^+$.
We set 
\[ I(x,y,t)= -  \tanh \left(\sigma \exp \left(-\frac t c \right) \right) b(\lambda,x,y) ,\]
where 
\begin{displaymath}
\begin{array}{c}
 b(\lambda,x_1,x_2) = \int_{-1}^{1} \int_{-1}^{1} K(x_1,x_2,y_1,y_2) dy_1 dy_2=\\
    =  \frac{\pi}{4 \lambda}  \left( \makebox{Erf} (\sqrt{\lambda} (1-x_1))+ \makebox{Erf} (\sqrt{\lambda} (1+x_1))\right) 
\left( \makebox{Erf} (\sqrt{\lambda} (1-x_2))+ \makebox{Erf} (\sqrt{\lambda} (1+x_2))\right) ,		
\end{array}
\end{displaymath}
where $Erf$ represents the Gaussian error function.

In this case, it is easy to check that the exact solution is 
\[V (\bar{x}, t)= \exp(- \frac t c) .\]

The initial condition is $V_0(\bar{x}) \equiv 1$.

For the space discretisation, we have used $k=4$, that is, 4 Gaussian nodes in each subinterval.
Since the discretisation error in space must be $O(h^8)$, we consider it negligible, compared
with the discretisation error in time.

With the following tests, we want to check that the discretisation error in time satisfies
the condition
\[ e_i= \|V_i -U_i \| = O (h_t^2).  \]  

The results are displayed in Table~1. We have used two different time steps, $h_t=0.02$
and $h_t=0.01$, and we have approximated the solution over the time interval $[0,0.1]$.
For the space discretisation, we have considered $N=24$, $m=12$. The equation parameters are 
$\lambda = \sigma=c=1$.

The discretisation errors $ e_i(h_t)=\|V_i -U_i \|$ are displayed at different moments $t_i$,
for different stepsizes $h_t$. We also present the ratios $e_i(2h_t)/e_i(h_t)$, which allow us
check the convergence order.
The ratios are  close to 4, which confirms the second order convergence. 
%%%%%%%%%%%%%%%%%%%%%%%%%%%%%%%%%%%%%%%%%%%%%%
\begin{table}
\begin{displaymath}
\begin{array} {|c|c|c|c|}
\hline
t  &   e_i(0.01) & e_i(0.02)  & e_i(0.02)/e_i(0.01)\\
\hline
0.02 & 6.66E-5  &   &    \\
0.03 & 7.24E-5  &   &    \\
0.04 & 7.46E-5  & 2.66E-4&  3.57  \\
0.05 & 7.56E-5  &    & \\
0.06 & 7.61E-5  & 2.91E-4 &   3.82  \\
0.07 & 7.65E-5  &         &       \\
0.08 & 7.69 E-5 & 3.01E-4 &  3.91  \\
0.09 & 7.72 E-5 &       &  \\
0.10 &  7.76E-5 &  3.06 E-4 &  3.94 \\    
\hline
\end{array}
\end{displaymath}
\caption{Numerical results for Example~1} 
\end{table}
%%%%%%%%%%%%%%%%%%%%%%%%%%%%%%%%%%%%%%%%%%%%%%%%%%%%%%%%%%

Now, in order to check the convergence of the space discretisation, we choose an example, where 
the time discretisation is exact (does not produce any error).

{ \bf Example 2.} In this example, the functions $K$ and $S$ are the same as in example 1. 
We set 
\[ I(x,y,t)= c+t-  \tanh (\sigma  t )  b(\lambda,x,y). \]
As in Example~1, $c=1$.
In this case, it is easy to check that the exact solution is 
\[V (\bar{x}, t)= t .\]

The initial condition is $V_0(\bar{x}) \equiv 0$. The difference operator \eqref{23} is exact
for linear functions of $t$, and this is why the scheme in this case does not have 
discretisation error in time. Therefore, the observed errors result from the space discretisation.
In this case, we are only considering the norm of the error at $t=0.1$.
The derivatives of $S$ and $K$ with respect to the space variables depend strongly
from the values of $\lambda$ and $\sigma$, and therefore these values should
influence the error of the space discretisation.  To check this, we consider 3 different cases:
$\lambda=1,\sigma=1$ ; $\lambda=1,\sigma=5$; and $\lambda=5, \sigma=5$, which are described
in tables 2,3 and 4, respectively.

Since we are using 4 Gaussian points in each subinterval, we expect that the error of
the space discretisation is $O(h^8)$ . Therefore, when we duplicate the number $N$ of gridpoints, the error
should decrease by a factor of approximately $2^8=256$.

In order to check the influence of interpolation error, for each $N$, we consider a set of 
 different values of $m$ (interpolation polynomial degree).

When $m$ increases from $12$ to $24$, the difference in accuracy is not significant.
This means that for values of $N$ up to $96$  it is enough to consider $m=12$.
When $\lambda$ or $\sigma$ increase we observe that the errors (for the same discretisation step)
also increase. This could be expected, since the discretisation error in space depends 
on the derivatives  $ \frac{ \partial ^i K(\bar{x},y_1,y_2)}
{\partial^i y_j}$ and $ \frac{ \partial ^i S(V)}
{\partial^i V}$, which increase with  $\lambda$ and $\sigma$, respectively.

%%%%%%%%%%%%%%%%%%%%%%%%%%%%%%      table 2
\begin{table}
\begin{displaymath}
\begin{array} {|c|c|c|c|c|c|}
\hline
m &   N=12 &  N=24  & e_{12}/e_{24} & N= 48 &  e_{24}/e_{48}\\
\hline
12 & 3.11 E-10   &   1.11E-12  &   280   & 3.997 E-15 &  278  \\
24 &             &   1.03E-12  &     &  4.413 E-15 &   234 \\
\hline
\end{array}
\end{displaymath}
\caption{Numerical results for Example~2, with $\lambda=1,\sigma=1$.  } 
\end{table}
%%%%%%
%%%%%%%%%%%%%%%%%%%%%%%%%%%%%%      table 3
\begin{table}
\begin{displaymath}
\begin{array} {|c|c|c|c|c|c|}
\hline
m &  N=24  &  N= 48 &  e_{24}/e_{48} &  N= 96 &  e_{48}/e_{96}  \\
\hline
12 & 1.62E-10   &   5.52E-13  &  293   &  2.36E-15  &  234  \\
24 & 1.69E-10   &   5.33E-13  &  317   &   2.22E-15 &  240  \\
\hline
\end{array}
\end{displaymath}
\caption{Numerical results for Example~2, with $\lambda=5,\sigma=1$.  } 
\end{table}
%%%%%%%%%%%%%%%%%%%%%%%%%%%%%%      table 4
\begin{table}
\begin{displaymath}
\begin{array} {|c|c|c|c|c|c|}
\hline
m &  N=24  &  N= 48 &  e_{24}/e_{48} &  N= 96 &  e_{48}/e_{96}  \\
\hline
12 & 7.31E-10   &   2.48E-12  &  295   &  9.38 E-15  & 264   \\
24 & 7.65E-10   &   2.40E-12  &  319   &   8.94E-15 &  268   \\
\hline
\end{array}
\end{displaymath}
\caption{Numerical results for Example~2, 
with $\lambda=5,\sigma=5$.  } 
\end{table}
%%%%%%%%%%%%%%%%%%%%%%%%%%%%%%%
{ \bf Example 3.} Finally let us consider an example where the potential distribution is not
constant, nor in time, neither in space. 
In this example,  the function $K$ has the same form as in the previous ones, but the forcing
function is
\[ I(x_1,x_2,t)= -\exp\left(-\frac{t}{c} \right) \beta(\lambda,\mu,x_1,x_2), \]
where
\[ \beta( \lambda, \mu, x_1,x_2)  = \int_0^1 \int_0^1  \exp \left(-\lambda ((x_1-y_1)^2 + (x_2-y_2)^2) -
\mu(y_1^2 + y_2^2)    \right) dy_1 dy_2 . \]
Let us consider a linear firing rate function $S(x)=x$.
If the we set the initial condition 
\[ V_0(x_1,x_2)= \exp \left(-\mu (x_1^2 + x_2^2)\right),\]
we conclude that the exact solution of this problem is
\[ V(x_1,x_2,t)= \exp \left(-\frac{t}{c} \right) \exp \left(-\mu (x_1^2 + x_2^2)\right).\]
We have computed the numerical solution by our method , over the time interval $[0,0.05]$,
with stepsize $h_t=0.01,0.005$ and $h_t=0.0025$. The parameters of the space discretisation are
$m=12$,$N=24$, therefore the error resulting from the space discretisation is in this case negligible,
 compared with the global error.
The error analysis for this example is displayed in table~5.
We see again that the convergence rate is in agreement with the theoretical results. 
%%%%%%%%%%%%%%%%%%%%%%%%%%%%%%      table 5
\begin{table}
\begin{displaymath}
\begin{array} {|c|c|c|}
\hline
h_t &  \|e_{h_t} \|_{\infty} &  \|e_{h_t} \|_{\infty}/\|e_{h_t/2} \|_{\infty} \\
\hline
0.01 & 7.66E-5   &   3.97  \\
0.005 & 1.93E-5   &   3.99 \\
0.0025 & 4.83E-6  &    \\
\hline
\end{array}
\end{displaymath}
\caption{Error norms and convergence rates for Example~3, 
with $\lambda=1,\sigma=1$.} 
\end{table}

\subsection{Neural Field Equation with Delay}
{ \bf Example 4.} In order to analyse the effect of delay, we now consider
equation \eqref{2} with the same data as in Example~3, but with some finite  propagation speed $v$.
We consider the initial condition $V_0(\bar{x}) = \exp \left(-\mu (x_1^2 + x_2^2)\right)$,
$\forall \bar{x} \in \Omega, t \in [-\tau_{max},0].$

In Fig. 1 some graphs of the solution are displayed. On the left-hand side, one can see the plots
of the solution at $t=0.5$, $t=1$, $t=1.5$, and $t=2.0$, for the non-delay case. On the right-hand side the corresponding plots are depicted, but for the delayed equation, when the propagation speed is $v=1$. The values of the remaining parameters
are $c=1$, $\mu=1$, $\lambda=1$. In both cases the stepsize 
in time is $h_t=0.1$, and the parameters of the space discretisation are $m=12$, $N=24$.
From this figure it is clear that as an effect of the delay, the decay of the solution in the 
case of the delayed equation is much slower.
\begin{figure}[p]
\includegraphics[clip=true, trim= 1cm 6cm 1cm 6cm,  width=14cm]{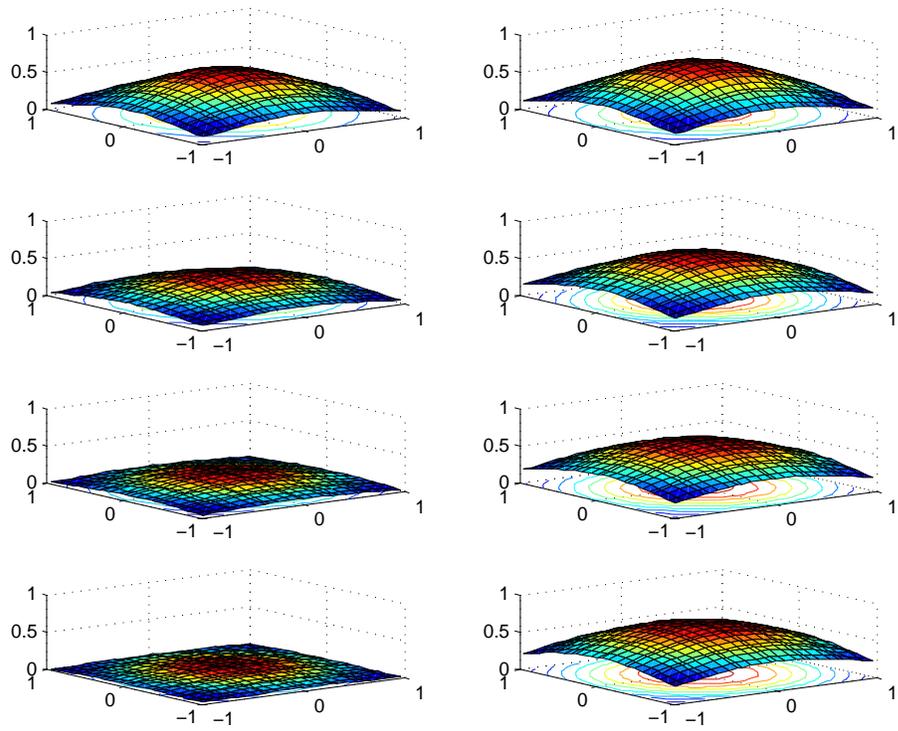} 
\caption{Plots of the solution  without delay (left) and with delay (right)}
\end{figure}

\section{Conclusions}
We have described and analysed a new numerical algorithm for computing approximate solutions of the two-dimensional neural field equations with delay. The stability, convergence and complexity of this algorithm have been analysed and a set of numerical examples has been presented. The numerical experiments are in agreement with the theoretical results and confirm that this algorithm can be successfully used for the solution of problems in Neuroscience and Robotics. 
 The main advantages of the described algorithm are its stability and accuracy, which is illustrated by the presented examples.
 Moreover, due to the use of a rank reduction technique, it is efficient when dealing with the two-dimensional case.

\section{Acknowledgements}
This research was supported by a Marie Curie Intra European Fellowship within the 7th European Community Framework Programme (PIEF-GA-2013-629496).

\end{document}